\documentclass[reqno]{alea2}
\usepackage{natbib}
\usepackage{fancyhdr}
\usepackage{graphicx}

\usepackage{amssymb}
\usepackage{amsmath}

\renewcommand{\cite}{\citet}

\makeatletter \@addtoreset{equation}{section} \makeatother

\renewcommand\theequation{\thesection.\arabic{equation}}
\renewcommand\thefigure{\thesection.\@arabic\c@figure}
\renewcommand\thetable{\thesection.\@arabic\c@table}

\theoremstyle{plain}
\newtheorem{theorem}{Theorem}[section]                                          \newtheorem{proposition}[theorem]{Proposition}                          
\newtheorem{lemma}[theorem]{Lemma}

\theoremstyle{definition}

\theoremstyle{remark}

\newcommand{\bb}[1]{{\mathbb #1}}

\begin{document}

%\date{, accepted }

\iffalse
\documentclass[reqno,11pt]{ams}

\usepackage{amsmath,amsthm}
\usepackage{amsfonts,amssymb}
\setlength{\textwidth}{157.5mm}
\setlength{\textheight}{24.5cm}
\setlength{\topmargin}{-1cm}
\setlength{\oddsidemargin}{0cm}
\setlength{\evensidemargin}{0cm}
\newtheorem{theorem}{Theorem}
\newtheorem{lemma}[theorem]{Lemma}
\newtheorem{proposition}[theorem]{Proposition}
\newtheorem{corollary}[theorem]{Corollary}
\newtheorem{example}{Example}
\newtheorem{definition}{Definition}
\newtheorem{conjecture}{Conjecture}
\fi

\def \diag{\hspace{1pt} {\rm diag}}
\def\P{{
\mathcal P}}
\def\Q{{\mathcal Q}}
\def \sh{{\rm sh \:}}
\def \sgn{{\rm sgn }}
\def \tr{{\rm tr }}
\def \pr{{\mathbb P}}
\def \Pm{{\bf P}}
\def \meanm{{\bf E}}
\def \mean{{\mathbb E}}
\def \Par{W}
\def \sgn{{\rm sgn }}
\def \N {{\bb N}}
\def \Z {{\bb Z}}
\def \S {{\mathcal S}}
\def \id {{\rm id}}
\def \R {{\bb R}}
\def \sp{{\rm sp}}
\def \Mat{{\mathbf H}}

\title[The largest-eigenvalue process for Wishart random matrices]%
{On the largest-eigenvalue process \\
for generalized Wishart random matrices}
\author{A. B. Dieker}
\address{Georgia Institute of Technology, Atlanta GA 30332, USA}
\email{ton.dieker@isye.gatech.edu}
\author{J. Warren}
\address{University of Warwick, Department of Statistics, Coventry, CV4 7AL, United Kingdom}
\email{j.warren@warwick.ac.uk}

\begin{abstract}
Using a change-of-measure argument, we prove an equality in law between the process of largest eigenvalues in a generalized Wishart random-matrix process and a last-passage percolation process.
This equality in law was conjectured by \cite{borodinpeche2008}.
\end{abstract}

\maketitle

%\vspace{-6mm} 
\section{Introduction}
The past decade has witnessed a surge of interest in connections between random matrices on the one hand
and applications to growth models, queueing systems, and last-passage percolation models on the other 
hand; standard references are \cite{baryshnikov:guequeues2001} and \cite{johansson:shapefluc2000}.
In this note we prove a result of this kind:
an equality in law between a process of largest eigenvalues for a family of Wishart random matrices
and a process of directed last-passage percolation times.

To formulate the main result, we  
construct two infinite arrays of random variables on an underlying measurable space,
along with a family $\{P^{\pi,\hat\pi}\}$ of probability measures parametrized by a 
positive $N$-vector $\pi$ and
a nonnegative sequence $\{\hat\pi_n: n\ge 1\}$.
The elements of the first array $\{A_{ij}: 1\le i\le N, j\ge 1\}$
are independent and $A_{ij}$ has a complex zero-mean Gaussian distribution with variance $1/(\pi_i+\hat\pi_j)$
under $P^{\pi,\hat\pi}$.
That is, both the real and complex part of $A_{ij}$ have zero mean and variance $1/(2\pi_i+2\hat\pi_j)$.
Write $A(n)$ for the $N\times n$ matrix formed by the first $n$ columns of $A$,
and define the matrix-valued stochastic process $\{M(n):n\ge 0\}$
%through the family $\{A(n): n\ge 0\}$ of `nested' matrices
by setting $M(n)=A(n)A(n)^*$ for $n\ge 1$ and by
letting $M(0)$ be the $N\times N$ zero matrix.
We call $\{M(n):n\ge 0\}$ a generalized Wishart random-matrix process, since
the marginals have a Wishart distribution if $\pi$ and $\hat\pi$ are identically one and zero, respectively.

The elements of the second array $\{W_{ij}: 1\le i\le N, j\ge 1\}$ are independent and 
$W_{ij}$ is exponentially distributed with parameter $\pi_i +\hat\pi_j$
under $P^{\pi,\hat\pi}$. We define
\[
Y(N,n) = \max_{P\in \Pi(N,n)} \sum_{(ij)\in P} W_{ij},
\]
where $\Pi(N,n)$ is the set of up-right paths from $(1,1)$ to $(N,n)$. The quantity
$Y(N,n)$ arises in last-passage percolation models as well as in series Jackson networks in queueing 
theory, see for instance \cite{diekerwarren:determinant2008} or \cite{johansson:meixner2007}.

The following theorem, a process-level equality in law 
between the largest eigenvalue of $M(n)$ and $Y(N,n)$, is the main result of this note.
Given a matrix $C$, we write $\sp(C)$ for its vector of eigenvalues, ordered decreasingly.

\begin{theorem}
\label{conj}
For any strictly positive vector $\pi$ and any nonnegative sequence $\hat\pi$,
the processes $\{\sp(M(n))_1 :n\ge 1\}$ and $\{Y(N,n):n\ge 1\}$ have the same distribution under
$P^{\pi,\hat\pi}$.
\end{theorem}

It is known from \cite{defosseux2008,forresterrains:jacobians2006} that this holds in the `standard' case, i.e., under the measure $P:=P^{(1,\ldots,1),(0,0,\ldots)}$.
In its stated generality, the theorem was conjectured by
\cite{borodinpeche2008}, who prove that the laws of $Y(N,n)$ and 
the largest eigenvalue of $M(n)$ coincide for fixed $n\ge 1$.
Our proof is based on a change-of-measure argument, which is potentially
useful to prove related equalities in law.

Throughout, we use the following notation. 
We let $\Mat_{N,N}$ be the space of all $N\times N$ Hermitian matrices, and
$\Par^N$ the set $\{x\in\R^N:x_1\ge \ldots\ge x_N\}$.
For $x,x'\in \Par^N$, we write $x\prec x'$ to mean that $x$ and $x'$ interlace
in the sense that
\[
x'_1\ge x_1\ge x'_2\ge x_2\ge \ldots\ge x_N'\ge x_N.
\]

\section{Preliminaries}
This section provides some background on
generalized Wishart random matrices, and introduces a Markov chain which plays
an important role in the proof of Theorem~\ref{conj}.

\subsection{The generalized Wishart random-matrix process}
\label{sec:wishart}
Under $P^{\pi,\hat\pi}$, the generalized Wishart process $\{M(n):n\ge 0\}$ from the introduction has independent increments since, for $m\ge 1$,
\begin{equation}
\label{eq:constrM}
M(m)= M(m-1) + \left(A_{im} \bar A_{jm}\right)_{1\le i,j\le N},
\end{equation}
where $\bar A_{jm}$ is the complex conjugate of $A_{jm}$.  
In particular, the matrix-valued increment has unit rank. 
The matrix $M(m)-M(m-1)$ can be parameterized by its diagonal elements together with the complex
arguments of $A_{im}$ for $1\le i\le N$; under $P^{\pi,\hat\pi}$, these are independent
and the former have exponential distributions
while the latter have uniform distributions on $[0,2\pi]$. 
(This fact is widely used in the Box-Muller method for computer generation of
random variables with a normal distribution.) 
Since the $i$-th diagonal element has an exponential distribution under $P^{\pi,\hat\pi}$ with
parameter $\pi_i+\hat\pi_m$, we obtain the following proposition.

\begin{proposition}
\label{prop:wishart}
For any $m\ge 1$, the $P^{\pi,\hat\pi}$-law of $M(m)-M(m-1)$ is absolutely continuous 
with respect to the $P$-law of $M(m)-M(m-1)$, and the Radon-Nikodym derivative is
\[
\prod_{i=1}^N(\pi_i+\hat\pi_m) \exp\left( -\sum_{i=1}^N (\pi_i+\hat\pi_m -1) (M_{ii}(m)-M_{ii}(m-1)) \right).
\]
\end{proposition}

\subsection{A Markov transition kernel}
\label{sec:markovchain}
We next introduce a time-inhomogeneous Markov transition kernel on $W^N$.
We shall prove in Section~\ref{sec:proof} that this kernel
describes the eigenvalue-process of the generalized Wishart random-matrix process
of the previous subsection.

In the standard case ($\pi \equiv 1$, $\hat\pi \equiv 0$), it follows 
from unitary invariance (see \cite[Sec.~5]{defosseux2008} or
\cite{forresterrains:jacobians2006}) that the process $\{\sp(M(n): n \geq 0\}$ is a 
homogeneous Markov chain.
Its one-step transition kernel $Q(z,\cdot)$ 
is the law of $\sp( \diag(z)+ G)$, where $G=\{g_i\bar{g}_j: 1\le i,j\le N\}$ is a rank one 
matrix determined by an $N$-vector $g$ of standard complex Gaussian random variables.
For $z$ in the interior of $W^N$, $Q(z,\cdot)$ is absolutely continuous with respect 
to Lebesgue measure on $W^N$ and can be written explicitly as in \cite[Prop.~4.8]{defosseux2008}:
\[
Q(z,dz') = \frac{\Delta(z')}{\Delta(z)} e^{-\sum_k (z_k'-z_k)} 1_{\{z\prec z'\}}dz',
\]
where $\Delta(z):=\prod_{1\le i<j\le N} (z_i-z_j)$ is the Vandermonde determinant.

\iffalse
Recall the Harish-Chandra-Itzykson-Zuber formula
\begin{equation}
\label{eq:HC}
 \int_U  \exp\left(-\tr(\diag(\pi) U\diag(z) U^*)\right) dU = c_{N}^{-1} \frac{\det\{e^{-\pi_iz_j}\}}{\Delta(\pi)\Delta(z)}, %\label{eq:RNdifficultterm}
\end{equation}
writing $dU$ for normalized Haar measure on the unitary group and $c_N$ for some (known) constant.
\fi
We use the Markov kernel $Q$ to define the aforementioned time-inhomogeneous Markov kernels, which
arise from the generalized Wishart random-matrix process.
For general $\pi$ and $\hat \pi$, we define the inhomogeneous transition 
probabilities $Q^{\pi,\hat\pi}_{n-1,n}$ via
\[
Q^{\pi,\hat\pi}_{n-1,n}(z,dz')=  
\prod_{i=1}^N (\pi_i+\hat\pi_n)  \frac{h_\pi(z^\prime)}{h_\pi(z)} e^{-(\hat\pi_n-1) 
\sum_{i=1}^N (z'_i-z_i)}   Q(z,dz^\prime),
\]
with
\begin{equation}
\label{eq:defh}
h_\pi(z) = \frac{\det\{e^{-\pi_iz_j}\}}{\Delta(\pi)\Delta(z)}.
\end{equation}
Note that $h_\pi(z)$ extends to a continuous function on $(0,\infty)^N \times W^N$ 
(this can immediately be seen as a consequence of the Harish-Chandra-Itzykson-Zuber formula,
see (\ref{eq:HC}) below).

One can verify that the $Q^{\pi,\hat\pi}$ are true Markov kernels by writing $1_{\{z\prec z'\}}=\det \{1_{\{z_i<z_j'\}}\}$ and applying the Cauchy-Binet formula
\[
\int_{W^N} \det \bigl\{ \xi_i(z_j)\bigr\} \det \bigl\{ \psi_j(z_i)\bigr\} dz = 
\det \left\{\int_{\R} \xi_i(z)\psi_j(z)dz\right\}.
\]

\section{The generalized Wishart eigenvalue-process}
\label{sec:proof}
In this section, we determine the law of the eigenvalue-process of generalized Wishart random-matrix process.
Although it is not essential to the proof of Theorem~\ref{conj},
we formulate our results in a setting where $\sp(M(0))$ is allowed to be nonzero.

Write $m_\mu$ for the `uniform distribution' on the set $\{M\in\Mat_{N,N}: \sp(M)=\mu\}$.
That is, $m_\mu$ is the unique probability 
measure invariant under conjugation by unitary matrices, or equivalently
$m_\mu$ is the law of $U\diag(\mu)U^*$ where $U$ is unitary and distributed according to (normalized) 
Haar measure. We define measures $P^{\pi,\hat\pi}_\mu$ by 
letting the $P^{\pi,\hat\pi}_\mu$-law of $\{M(n)-M(0): n\ge 0\}$ be equal to
the $P^{\pi,\hat\pi}$-law of $\{M(n):n\ge 0\}$,
and letting the $P^{\pi,\hat\pi}_\mu$-distribution of $M(0)$ be 
independent of $\{M(n)-M(0):n\ge 0\}$
and absolutely continuous with respect to $m_\mu$ with Radon-Nikodym derivative 
\begin{equation}
\label{eq:RNforM0}
%(2\pi)^{N/2} \prod_{j=1}^{N-1} [j!] 
\frac{c_N}{h_\pi(\mu)}e^{-\sum_{i=1}^N \mu_i} 
%\frac{\Delta(\pi)\Delta(\mu)}{\det\{e^{-\pi_i \mu_j} \}} 
\exp(-\tr[(\diag(\pi)-I) M(0)]),
\end{equation}
where $c_N$ is a constant depending only on the dimension $N$ and $I$ is the identity matrix.
Recall that $h_\pi(\mu)$ is defined in (\ref{eq:defh}).
That this defines the density of a probability measure for all $\pi$ and
$\mu$ follows immediately from the Harish-Chandra-Itzykson-Zuber formula 
(e.g., \cite[App.~A.5]{MR2129906})
\begin{equation}
\label{eq:HC}
\int_U  \exp\left(-\tr(\diag(\pi) U\diag(\mu) U^*)\right) dU = c_{N}^{-1} h_\pi(\mu),
\end{equation}
writing $dU$ for normalized Haar measure on the unitary group. 
Throughout, we abbreviate $P^{(1,\ldots,1),(0,0,\ldots)}_\mu$ by $P_\mu$.
Note that the $P^{\pi,\hat\pi}_\mu$-law and the $P^{\pi,\hat\pi}$-law of $\{M(n):n\ge 0\}$ coincide
if $\mu=0$.

%Consider a distribution under which the infinite array $A$ has the same distribution as under $P^{\pi,\hat\pi}$,
%and under which $M(0)$ is independent of $A$ with the same distribution as $M(N)$ ;
%the process $\{M(n):n\ge 0\}$ is constructed from $M(0)$ and $A$ as in (\ref{eq:constrM}).
%We let $P^{\pi,\hat\pi}_\mu$ be this distribution conditioned on $\sp(M(0))=\mu$.

The following theorem specifies the $P^{\pi,\hat\pi}_\mu$-law of $\{\sp(M(n)):n\ge 0\}$.
% is the
%same as the $P^{\pi,\hat\pi}$-law of $\{Z(n):n\ge 0\}$, and 
%in conjuction with Lemma~\ref{lem:shape}, it implies Theorem~\ref{conj}.

\begin{theorem}
\label{thm:main}
For any $\mu\in W^N$,
$\{\sp(M(n)): n\ge 0\}$ is an inhomogeneous Markov chain on $W^N$ under $P^{\pi,\hat\pi}_\mu$,
and it has the $Q^{\pi,\hat\pi}_{n-1,n}$ of Section~\ref{sec:markovchain}
for its one-step transition kernels.
\end{theorem}
\proof
Fix some $\mu\in W^N$. The key ingredient in the proof is a change of measure argument.
We know from \cite{defosseux2008} or \cite{forresterrains:jacobians2006}
%and Baryshnikov~\cite{baryshnikov:guequeues2001} 
that Theorem~\ref{thm:main} holds for the `standard' case $\pi=(1,\ldots,1)$, $\hat\pi\equiv 0$.
%From (\ref{eq:equationQn}),
%it is readily seen that the function $f^\nu$ on $\Mat_{N+1,N+1}$ defined by
%\[
%f^\nu(M) = \exp\left(-\tr(\diag(\nu) M) + \tr( M)\right)
%\]
%is an eigenfunction for $Q_n$ (with an eigenvalue depending on $\nu$).

\newcommand{\MN}{M}
Writing $P^{\pi,\hat\pi}_n$ and $P_n$ for the distribution of
$(\MN(0),\ldots,\MN(n))$ under $P^{\pi,\hat\pi}_\mu$ and $P_\mu$ respectively,
we obtain from Section~\ref{sec:wishart} that for $n\ge 0$,
\begin{eqnarray*}
%\label{eq:RNderivativemat}
\lefteqn{\frac{dP_n^{\pi,\hat\pi}}{dP_n} (\MN(0),\ldots,\MN(n))} \\&=&
C_{\pi,\hat\pi}(n,N) \frac{c_N}{h_\pi(\mu)} 
e^{-\sum_{i=1}^N \mu_i} \\
&&\mbox{}\, \times\, \exp\left(-\tr((\diag(\pi)-I)\MN(n))
-\sum_{m=1}^{n}\hat\pi_m \tr(\MN(m)-\MN(m-1))\right),%\nonumber
\end{eqnarray*} 
where $C_{\pi,\hat\pi}(n,N) = \prod_{i=1}^N \prod_{j=1}^{n} (\pi_i+\hat\pi_j)$.
Let the measure $p_n^{\pi,\hat\pi}$ (and $p_n$) be the restriction of $P_n^{\pi,\hat\pi}$ (and $P_n$) to the $\sigma$-field generated by $(\sp(\MN(0)),\ldots,\sp(\MN(n)))$. Then we obtain for $n\ge 0$,
\begin{eqnarray*}
\lefteqn{\frac{dp_n^{\pi,\hat\pi}}{dp_n} (\sp(\MN(0)),\ldots,\sp(\MN(n)))} \\&=& 
\mean_{P_\mu}\left[\left.\frac{dP_n^{\pi,\hat\pi}}{dP_n} (\MN(0),\ldots,\MN(n))\right| \sp(\MN(0)),\ldots,
\sp(\MN(n))\right],
\end{eqnarray*}
where $\mean_{P_\mu}$ denotes the expectation operator with respect to $P_\mu$.
Since the $P_\mu$-distribution of $(\MN(0),\ldots,\MN(n))$ given the spectra is invariant under
componentwise conjugation by a unitary matrix $U$, we have 
for $\mu\equiv\mu^{(0)}\prec\mu^{(1)}\prec\ldots \prec \mu^{(n)}$,
\begin{eqnarray*}
\lefteqn{\mean_{P_\mu}\left[ \exp\left(-\tr(\diag(\pi) \MN(n))\right)\left| \sp(\MN(0))=\mu^{(0)},\ldots,
\sp(\MN(n))=\mu^{(n)}\right. \right]} \nonumber \\
&=& \int_U  \exp\left(-\tr(\diag(\pi) U\diag(\mu^{(n)}) U^*)\right) dU \nonumber
\hspace{40mm}
\\&=& c_{N}^{-1} h_\pi(\mu^{(n)}), %\label{eq:RNdifficultterm}
\end{eqnarray*}
where the second equality is the Harish-Chandra-Itzykson-Zuber
formula.
From the preceding three displays in conjunction with $\tr(M)=\sum_i \sp(M)_i$, we conclude that
\begin{eqnarray*}
\lefteqn{\frac{dp^{\pi,\hat\pi}_n}{dp_n} (\mu,\mu^{(1)},\ldots,\mu^{(n)})}\\ &=& 
C_{\pi,\hat\pi}(n,N) 
\frac{h_\pi(\mu^{(n)})}{h_\pi(\mu)}
%\\&&\mbox{}\times
\exp\left(-\sum_{i=1}^N \sum_{r=1}^n \hat\pi_r \left[\mu_i^{(r)}-\mu_i^{(r-1)}\right]+
\sum_{i=1}^N [\mu_i^{(n)}-\mu_i] \right).
\end{eqnarray*}
Since $\sp(M(\cdot))$ is a Markov chain with transition kernel $Q$ under $P_\mu$, we have
\begin{eqnarray*}
\lefteqn{P^{\pi,\hat\pi}_\mu(\sp(\MN(1))\in d\mu^{(1)}, \ldots, \sp(\MN(n))\in d\mu^{(n)} )} 
\\&=&\frac{dp_n^{\pi,\hat\pi}}{dp_n} (\mu,\mu^{(1)},\ldots,\mu^{(n)}) 
P_\mu(\sp(\MN(1))\in d\mu^{(1)}, \ldots, \sp(\MN(n))\in d\mu^{(n)} )\\ &=&
\frac{dp_n^{\pi,\hat\pi}}{dp_n} (\mu,\mu^{(1)},\ldots,\mu^{(n)}) Q(\mu,d\mu^{(1)})\cdots 
Q(\mu^{(n-1)},d\mu^{(n)})\\
&=& Q^{\pi,\hat\pi}_{0,1} (\mu, d\mu^{(1)}) 
Q^{\pi,\hat\pi}_{1,2} (\mu^{(1)}, d\mu^{(2)}) \cdots Q^{\pi,\hat\pi}_{n-1,n} (\mu^{(n-1)}, d\mu^{(n)}),
\end{eqnarray*}
%where the transition density $p(\cdot,\cdot)$ is defined in (\ref{eq:kernelstd}).
%Consequently, if $\mu\prec\mu^{(1)}$,
%\begin{eqnarray*}
%\lefteqn{P^{\pi,\hat\pi}_\mu \left(\sp(M(1))\in d\mu^{(1)}, \ldots, \sp(M(n))\in d\mu^{(n)} \right)} \\
%&=& P^{\pi,\hat\pi}\left(\sp(\MN(1))\in d\mu^{(1)}, \ldots, \sp(\MN(n))\in d\mu^{(n)} \left|\sp(\MN(0)) %=\mu\right.\right) \\
%&=& \frac{dp_n^{\pi,\hat\pi}}{dp_n} (\mu,\mu^{(1)},\ldots,\mu^{(n)}) p(\mu,\mu^{(1)})\cdots 
%p(\mu^{(n-1)},\mu^{(n)})d\mu^{(1)}\cdots d\mu^{(n)}\left/ \frac{dp_0^{\pi,\hat\pi}}{dp_0} 
%(\mu)\right. .
%\end{eqnarray*}
%We have thus shown that the $P^{\pi,\hat\pi}_\mu$-law of $(\sp(M(1)),\ldots,\sp(M(n)))$ is absolutely
%continuous with respect to Lebesgue measure, and that it has density
%\[
%\prod_{i=1}^N\prod_{j=1}^n (\pi_i+\hat\pi_j) 
%\frac{\det\{e^{-\pi_i\mu^{(n)}_j}\}}{\det\{e^{-\pi_i\mu_j}\}}
%\exp\left(-\sum_{i=1}^N \sum_{r=1}^n \hat\pi_r \left[\mu_i^{(r)}-\mu_i^{(r-1)}\right]\right)
%1_{\{\mu\prec\mu^{(1)}\prec\ldots\prec \mu^{(n)}\}}.
%\]
the last equality being a consequence of the definition of $Q^{\pi,\hat\pi}_{k-1,k}$
and the expression for $dp_n^{\pi,\hat\pi}/dp_n$.~\endproof

\section{Robinson-Schensted-Knuth and the proof of Theorem~\ref{conj}}
This section explains the connection between the infinite array $\{W_{ij}\}$ of the introduction and
the Markov kernels $Q^{\pi,\hat\pi}_{n-1,n}$.
In conjunction with Theorem~\ref{thm:main},
these connections allow us to prove Theorem~\ref{conj}.

\subsection*{The RSK algorithm}
The results in this section rely on a combinatorial
mechanism known as the Robinson-Schensted-Knuth (RSK) algorithm.
This algorithm generates from a $p \times q$ matrix with nonnegative entries
a triangular array $\mathbf x = \{x_i^j: 1\le j\le p, 1\le i \le j\}$ called a Gelfand-Tsetlin (GT) pattern. 
A GT pattern with $p$ levels $x^1,\ldots,x^p$ is an array for which the coordinates satisfy the inequalities
\[
x^k_k\le x_{k-1}^{k-1}\leq x^k_{k-1}\leq x^{k-1}_{k-2}
\le \ldots\leq x_2^k\leq x_1^{k-1}\leq x_1^k
\]
for $k=2,\ldots,p$.
%We write $\mathbf{K}^N$ for the set of all $\mathbf x$ satisfying the above constraint
If the elements of the matrix are integers, then a GT pattern
can be identified with a so-called semistandard Young tableau, and the bottom row $x^p=\{x^p_i;
1 \le i \le p\}$ of the GT pattern corresponds to the shape of the Young tableau.
We write $\mathbf K_p$ for the space of all GT patterns $\mathbf x$ with $p$ levels.

By applying the RSK algorithm with row insertion to an infinite array
$\{\xi_{ij}: 1\le i\le N, 1\le j\le n\}$ for $n=1,2,\ldots$, we obtain 
a sequence of GT patterns $\mathbf x(1),\mathbf x(2),\ldots$.
It follows from properties of RSK that
\begin{equation}
\label{eq:X1Nmax}
x_1^N(n)  = \max_{P\in \Pi(N,n)} \sum_{(ij)\in P} \xi_{ij},
\end{equation}
where $\Pi(N,n)$ is the set of up-right paths from $(1,1)$ to $(N,n)$ as before.
Details can be found in, e.g.,~\cite{johansson:shapefluc2000} or 
\cite[case A]{diekerwarren:determinant2008}. 

Greene's theorem generalizes (\ref{eq:X1Nmax}), and gives similar expressions for each component
of the pattern $x_i^j(n)$, see for instance Chapter 3 of \cite{fulton:youngtableaux1997} or
Equation (16) in \cite{MR2053054}.
%To explain this, observe that the RSK algorithm applied to $\{\xi_{ij}: 1\le i\le k, 1\le j\le n\}$ 
%for some $k=1,\ldots,N$ yields the GT pattern $\{x_i^j(n): 1\le i\le k, 1\le j\le k\}$ consisting 
%of the top $k$ levels of $\mathbf x(n)$.
As a consequence of these, we can consider the RSK algorithm for real-valued $\xi_{ij}$
and each $\mathbf x(n)$ is then a continuous function of the input data.

We remark that the RSK algorithm can also be started from a given initial GT pattern $\mathbf x(0)$. 
If RSK is started from the null pattern, it reduces to the standard algorithm
and we set $\mathbf x(0)=0$.

\subsection*{The bijective property of RSK}
RSK has a bijective property which has important probabilistic consequences
for the sequence of GT patterns constructed from specially chosen 
random infinite arrays.
Indeed, suppose that $\{\xi_{ij}: 1\le i\le N, j\ge 1\}$ is a family of independent random variables
with $\xi_{ij}$ having a geometric distribution on $\Z_+$ with parameter $a_i b_j$, where 
$\{a_i:1\le i\le N\}$ and $\{b_j: j\ge 1\}$ are two sequences taking values in $(0,1]$.
Write $\{\mathbf X(n):n\ge 0\}$ for the sequence of GT patterns constructed from $\xi$.

Using the bijective property of RSK 
it can be verified that the bottom rows $\{X^N(n):n\ge 0\}$ of the GT patterns
evolve as an inhomogeneous Markov chain with transition probabilities
\begin{equation}
\label{eq:transprobdiscr}
P_{n-1,n}(x,x') = \prod_{i=1}^N(1-a_ib_n) \frac{s_{x'}(a)}{s_{x}(a)} b_n^{\sum_{i=1}^N (x_i'-x_i)}
1_{\{0\le x\prec x'\}},
\end{equation}
where $s_\lambda(a)$ is the Schur polynomial corresponding to a partition $\lambda$:
\[
s_\lambda(a) = \sum_{\mathbf x\in\mathbf K_N: x^N=\lambda} a^{\mathbf x},
\]
with the weight $a^{{\mathbf x}}$ of a GT pattern $\mathbf x$ being defined as
\[
a^{\mathbf{x}}= a_1^{x^1_1} \prod_{k=2}^N a_k^{ \sum
x^k_i- \sum x^{k-1}_i}.
\]
This is proved in \cite{oconnell:conditionedRSK2003} in the special case with $b_j=1$ for all $j$,
and the argument extends straightforwardly; see also \cite{forresternagao2008}.

Non-null initial GT patterns generally do {\em not} give rise to Markovian bottom-row processes.
Still, the inhomogeneous Markov chain of bottom rows
can be constructed starting from a given initial partition $\lambda$ with at most $N$ parts
by choosing $\mathbf X(0)$ suitably from the space of a GT patterns with bottom row $\lambda$:
$\mathbf X(0)$ should be independent of the family $\{\xi_{ij}\}$ with probability mass function
\[
p({\mathbf x}) = \frac{a^{{\mathbf x}}}{s_\lambda(a)}.
\] 

\subsection*{Exponentially distributed input data}
We now consider the sequence of GT patterns $\{\mathbf X_{L}(n): n\ge 0\}$
arising from setting $a_i = 1-\pi_i/L$ and $b_j=1-\hat \pi_j/L$ in the above setup, 
and we study the regime $L\to\infty$ after rescaling suitably.
In the regime $L\to\infty$, the input variables $\{\xi_{ij}/L\}$ (jointly) converge in distribution
to independent exponential random
variables, the variable corresponding to $\xi_{ij}/L$ having parameter $\pi_i+\hat\pi_j$.
Thus, the law of the input array $\xi$ converges weakly to the 
$P^{\pi,\hat\pi}$-law of the array $\{W_{ij}: 1\le i\le N, j\ge 1\}$ from the introduction.
Refer to \cite{MR2053054} and \cite{johansson:shapefluc2000} for related results on this regime.

By the aforementioned continuity of the RSK algorithm and the continuous-mapping theorem,
$\{\mathbf X_{L}(n)/L: n\ge 0\}$ converges in distribution to a process $\{\mathbf Z(n):n\ge 0\}$
taking values in GT patterns with $N$ levels.
As a consequence of the above results in a discrete-space setting, we get from (\ref{eq:X1Nmax}) that
\[
Z_1^N(n)  = \max_{P\in \Pi(N,n)} \sum_{(ij)\in P} W_{ij}.
\]
Moreover, the process of bottom rows $\{Z^N(n):n\ge 0\}$ is an inhomogeneous Markov chain for which
its transition mechanism can be found by letting $L\to\infty$ in 
(\ref{eq:transprobdiscr}): %This yields the following lemma.

\begin{lemma}
\label{lem:RSK}
Under $P^{\pi,\hat\pi}$, the process $\{Z^N(n):n\ge0\}$ is an inhomogeneous Markov chain on $W^N$,
and it has the $Q^{\pi,\hat\pi}_{n-1,n}$ of Section~\ref{sec:markovchain}
for its one-step transition kernels.
\end{lemma}

A similar result can be obtained given a non-null initial bottom row $\mu\in W^N$.
In case the components of $\mu$ are distinct, the distribution of the initial pattern $\mathbf Z(0)$ 
should then be absolutely continuous
with respect to Lebesgue measure on $\{\mathbf z\in \mathbf K_N: z^N=\mu\}$ with density
\[
\frac{\Delta(\pi)}{\det\{e^{-\pi_i \mu_j}\}} c^{\mathbf z},
\]
where $c=(e^{-\pi_1},\ldots,e^{-\pi_N})$.
\vspace{-2mm}

\subsection*{Proof of Theorem~\ref{conj}}
We now have all ingredients to prove Theorem~\ref{conj}.
We already noted that  $Z_1^N(n)$ equals $Y(N,n)$.
Thus, for any strictly positive vector $\pi$ and any nonnegative sequence $\hat\pi$, 
$\{Y(N,n):n\ge 1\}$ has the same $P^{\pi,\hat\pi}$-distribution 
as $\{Z_1(n): n\ge 1\}$.
In view of Theorem~\ref{thm:main} and Lemma~\ref{lem:RSK}, in turn this
has the same $P^{\pi,\hat \pi}$-distribution as
the largest-eigenvalue process $\{\sp(M(n))_1:n\ge 1\}$. This proves Theorem~\ref{conj}.

{\small
\bibliography{../../../bibdb}
\bibliographystyle{alea2}
}

\end{document}